\documentstyle[12pt, amsfonts]{amsart}
\begin{document}
\def\R{{\mathbb R}}
\def\Z{{\mathbb Z}}
\def\C{{\mathbb C}}
\newcommand{\trace}{\rm trace}
\newcommand{\Ex}{{\mathbb{E}}}
\newcommand{\Prob}{{\mathbb{P}}}
\newcommand{\E}{{\cal E}}
\newcommand{\F}{{\cal F}}
\newtheorem{df}{Definition}
\newtheorem{theorem}{Theorem}
\newtheorem{lemma}{Lemma}
\newtheorem{pr}{Proposition}
\newtheorem{co}{Corollary}
\def\n{\nu}
\def\sign{\mbox{ sign }}
\def\a{\alpha}
\def\N{{\mathbb N}}
\def\A{{\cal A}}
\def\L{{\cal L}}
\def\X{{\cal X}}
\def\F{{\cal F}}
\def\c{\bar{c}}
\def\v{\nu}
\def\d{\delta}
\def\diam{\mbox{\rm dim}}
\def\vol{\mbox{\rm Vol}}
\def\b{\beta}
\def\t{\theta}
\def\l{\lambda}
\def\e{\varepsilon}
\def\colon{{:}\;}
\def\pf{\noindent {\bf Proof :  \  }}
\def\endpf{ \begin{flushright}
$ \Box $ \\
\end{flushright}}

\title[Stability in volume comparison]{Stability in the Busemann-Petty and Shephard problems}

\author{Alexander Koldobsky}

\address{Department of Mathematics\\ 
University of Missouri\\
Columbia, MO 65211}

\email{koldobskiya@@missouri.edu}

\begin{abstract} A comparison problem for volumes of convex bodies asks whether inequalities
$f_K(\xi)\le f_L(\xi)$ for all $\xi\in S^{n-1}$ imply that $\vol_n(K)\le \vol_n(L),$
where $K,L$ are convex bodies in $\R^n,$ and $f_K$ is a certain geometric characteristic of $K.$
By linear stability in comparison problems we mean that there exists a constant $c$ such that  for
every $\e>0$, the inequalities $f_K(\xi)\le f_L(\xi)+\e$ for
all $\xi\in S^{n-1}$ imply that $\left(\vol_n(K)\right)^{\frac{n-1}n}\le \left(\vol_n(L)\right)^{\frac{n-1}n}+c\e.$ 

We prove such results in the settings
of the Busemann-Petty and Shephard problems and their generalizations. 
We consider the section function $f_K(\xi)=S_K(\xi)=\vol_{n-1}(K\cap \xi^\bot)$ and the projection function
$f_K(\xi)=P_K(\xi)=\vol_{n-1}(K\vert\xi^\bot),$ where $\xi^\perp$ is the central hyperplane perpendicular to $\xi,$
and $K\vert\xi^\bot$ is the orthogonal projection of $K$ to $\xi^\bot.$
In these two cases we prove linear stability under additional conditions that
$K$ is an intersection body or $L$ is a projection body, respectively. 
Then we consider other functions $f_K,$ which allows to remove the additional conditions on the bodies 
in higher dimensions.
\end{abstract}  
\maketitle

\section{Introduction}

A typical comparison problem for the volume of convex bodies
asks whether inequalities $$f_K(\xi)\le f_L(\xi), \qquad  \forall \xi\in S^{n-1}$$ imply $\vol_n(K)\le \vol_n(L)$ 
for any $K,L$ from a certain class of origin-symmteric convex bodies in $\R^n,$ where $f_K$ 
is a certain geometric characteristic of $K$ and $\vol_n$ is the $n$-dimensional volume. 

If $f_K=S_K$ is the section function of $K$ defined by
$$S_K(\xi)=\vol_{n-1}(K\cap\xi^\perp), \qquad \xi\in S^{n-1},$$
where $\xi^\perp $ is the central hyperplane in $\R^n$ orthogonal to $\xi$,
the corresponding comparison question is the matter of the Busemann-Petty 
problem, raised in 1956 in \cite{BP} and solved in the end of the 1990's as the result 
of a sequence of papers \cite{LR}, \cite{Ba}, \cite{Gi}, \cite{Bo}, 
\cite{L}, \cite{Pa}, \cite{G1}, \cite{G2}, \cite{Z1}, \cite{Z2}, \cite{K2}, \cite{K3}, \cite{Z3},
\cite{GKS} ; see \cite[p. 3]{K4} or \cite[p. 343]{G3} for the history of the solution.
The answer is affirmative if $n\le 4$, and it is negative if $n\ge 5.$

Another example is the Shephard problem with $f_K=P_K$ being the projection function
$$P_K(\xi)=\vol_{n-1}(K\vert\xi^\perp), \qquad \xi\in S^{n-1},$$
where $K\vert\xi^\perp$ is the orthogonal projection of $K$ to the hyperplane $\xi^\perp.$
The Shephard problem was posed in 1964 in \cite{Sh} and solved soon after that by Petty \cite{Pe} and 
Schneider \cite{S1}. The answer if affirmative only  in dimension 2. 

Since the answers to the Busemann-Petty and Shephard problems are  negative in most dimensions,  
one may ask what information about the functions $S_K$ and $P_K$ does allow to compare the volumes 
in all dimensions. In the section case an answer to this question was given in \cite{KYY}:  
for two origin-symmetric infinitely smooth bodies $K,L$ in $\R^n$ and $\alpha\in [n-4,n-1)$ 
the inequalities
\begin{eqnarray}\label{eqn:condition1}
(-\Delta)^{\alpha/2} S_K(\xi)\le (-\Delta)^{\alpha/2} S_L(\xi), \qquad \forall \xi\in S^{n-1}
\end{eqnarray}
imply that $\vol_n(K)\le \vol_n(L),$ while for $\alpha<n-4$ this is not necessarily
true.
Here $\Delta$ is the Laplace operator on $\R^n$, and the fractional powers of the
Laplacian are defined by
$$
(-\Delta)^{\alpha/2}f = \frac{1}{(2\pi)^n}( |x|_2^\alpha  \hat{f}(x))^\wedge,
$$
where  the Fourier transform is considered in the sense of distributions, $|x|_2$ stands for
the Euclidean norm in $\R^n$, and the functions $S_K$ and $S_L$ are extended in (\ref{eqn:condition1}) 
to homogeneous functions of degree -1 on the whole $\R^n.$ This result contains the solution to 
the original Busemann-Petty problem as a particular case and means that one has to differentiate 
the section functions at least $n-4$ times in order to compare the $n$-dimensional volumes. 

The situation is different for projections where a similar extension does not directly generalize
the solution to Shephard's problem. Yaskin \cite{Y} proved that for $\alpha \in [n,n+1)$ the inequalities
\begin{eqnarray}\label{eqn:condition2}
(-\Delta)^{\alpha/2} P_K(\xi)\ge (-\Delta)^{\alpha/2} P_L(\xi), \qquad \forall \xi\in S^{n-1}
\end{eqnarray}
imply that $\vol_n(K) \le \vol_n(L),$ where the projection functions are extended to homogeneous functions 
of degree 1 on the whole $\R^n.$ The latter result is no longer true for $\alpha\in [n-2,n),$ which
would be a natural extension of the solution to the original Shephard's problem.
\medbreak
In this article we study the flexibility of the results mentioned above.  By {\it linear stability} in a comparison result
we mean that there exists a constant $c$ such that for any $K, L$ from certain classes of convex bodies
and every $\e>0$ the inequalities 
$$f_K(\xi)\le f_L(\xi) + \e, \qquad  \forall \xi\in S^{n-1}$$
imply 
$$\vol_n(K)^{\frac{n-1}n} \le \vol_n(L)^{\frac{n-1}n} + c\e.$$
We also consider {\it linear separation} in comparison problems, where we are looking for
a constant $c$ such that for any $K, L$ from certain classes of convex bodies
and every $\e>0$ the inequalities 
$$f_K(\xi)\le f_L(\xi) - \e, \qquad  \forall \xi\in S^{n-1}$$
imply 
$$\vol_n(K)^{\frac{n-1}n} \le \vol_n(L)^{\frac{n-1}n} - c\e.$$

We first prove linear stability and separation for the section function $f_K=S_K$ under the additional assumption
that $K$ is an intersection body.  In the stability result the constant $c=1,$ but in the case of separation $c$ depends on
the inradius of $K$ and on the dimension $n.$  Since every origin-symmetric convex body in $\R^n,\ 2\le n \le 4$ 
is an intersection body, in these dimensions the results apply to arbitrary origin-symmetric convex bodies 
$K,L.$  An easy consequence (see Corollary \ref{n4}) is that for $2\le n \le 4$ and any origin-symmetric convex bodies $K,L$ in $\R^n,$ 
$$\left|\vol_n(K)^{\frac{n-1}n} - \vol_n(L)^{\frac{n-1}n}\right| \le \|S_K-S_L\|_{C(S^{n-1})}.$$

We also prove linear stability and separation for the projection function $f_K=P_K$ under the additional assumption
that $L$ is a projection body.  Here in the stability result the constant $c$ depends on $n$ and on the circumradius of $L,$
while in the case of separation we have $c=\sqrt{1/e}.$

In order to remove the additional assumptions on the bodies and make the results work in general in higher
dimensions, we prove linear stability and separation in the results from \cite{KYY} and \cite{Y} mentioned above.
We consider the cases where  $$f_K = (-\Delta)^{\alpha/2} S_K,\qquad \alpha\in [n-4,n-1)$$
and $$f_K = (-\Delta)^{\alpha/2} P_K,\qquad \alpha\in [n,n+1),$$ and $K,L$ are arbitrary infinitely smooth convex 
bodies in $\R^n.$  In the stability case the constant $c$ for sections depends only on $\alpha$ and $n$, while 
for projections the constant also depends on the circumradius of $L.$ In the separation case, $c$ depends 
only on $\alpha$ and $n$ for projections, and also depends on the inradius of $K$ for sections.

In most cases we employ the techniques of the Fourier analytic  approach to sections and projections that has recently been developed;
see \cite{K4} and \cite{KY}. We use a more geometric Radon transform approach in the case $f_K=S_K$ to show
the variety of methods, it is also possible to solve this case with the Fourier transform.

\section{Stability and separation for sections}

We say that a closed bounded set $K$ in $\R^n$ is a {\it star body}  if 
every straight line passing through the origin crosses the boundary of $K$ 
at exactly two points different from the origin, the origin is an interior point of $K,$
and the {\it Minkowski functional} 
of $K$ defined by 
$$\|x\|_K = \min\{a\ge 0:\ x\in aK\}$$
is a continuous function on $\R^n.$ 

The {\it radial function} of a star body $K$ is defined by
$$\rho_K(x) = \|x\|_K^{-1}, \qquad x\in \R^n.$$
If $x\in S^{n-1}$ then $\rho_K(x)$ is the radius of $K$ in the
direction of $x.$

Writing the volume of $K$ in polar coordinates, one gets
\begin{equation} \label{polar-volume}
\vol_n(K)
=\frac{1}{n} \int_{S^{n-1}} \rho_K^n(\theta) d\theta=
\frac{1}{n} \int_{S^{n-1}} \|\theta\|_K^{-n} d\theta.
\end{equation}

The {\it spherical Radon transform} 
$R:C(S^{n-1})\mapsto C(S^{n-1})$  
is a linear operator defined by
$$Rf(\xi)=\int_{S^{n-1}\cap \xi^\bot} f(x)\ dx,\quad \xi\in S^{n-1}$$
for every function $f\in C(S^{n-1}).$

The polar formula for the volume of a hyperplane section expresses this volume 
in terms of the spherical Radon transform (see for example \cite[p.15]{K4}):
\begin{equation} \label{volume=spherradon}
S_K(\xi)= \vol_{n-1}(K\cap \xi^\bot) = 
\frac{1}{n-1} R(\|\cdot\|_K^{-n+1})(\xi).
\end{equation}

The spherical Radon 
transform is self-dual (see \cite[Lemma 1.3.3]{Gr} ):
for any functions $f,g\in C(S^{n-1})$
\begin{equation} \label{selfdual}
\int_{S^{n-1}} Rf(\xi)\ g(\xi)\ d\xi = \int_{S^{n-1}} f(\xi)\ Rg(\xi)\ d\xi.
\end{equation}

The spherical Radon transform can be extended to measures. 
Let $\mu$ be a finite Borel measure on $S^{n-1}.$
We define the spherical Radon transform of $\mu$ as a functional $R\mu$ on
the space $C(S^{n-1})$ acting by
$$(R\mu,f)= (\mu, Rf)= \int_{S^{n-1}} Rf(x) d\mu(x).$$
By Riesz's characterization of continuous linear functionals on the
space $C(S^{n-1})$, 
$R\mu$ is also a finite Borel measure on $S^{n-1}.$ If $\mu$ has 
continuous density $g,$ then by (\ref{selfdual}) the 
Radon transform of $\mu$ has density $Rg.$

The class of intersection bodies was introduced by Lutwak \cite{L}.
Let $K, L$ be origin-symmetric star bodies in $\R^n.$ We say that $K$ is the 
intersection body of $L$ if the radius of $K$ in every direction is 
equal to the $(n-1)$-dimensional volume of the section of $L$ by the central
hyperplane orthogonal to this direction, i.e. for every $\xi\in S^{n-1},$
\begin{equation} \label{intbodyofstar}
\rho_K(\xi)= \|\xi\|_K^{-1} = \vol_{n-1}(L\cap \xi^\bot).
\end{equation} \index{intersection body of a star body}
All the bodies $K$ that appear as intersection bodies of different star bodies
form {\it the class of intersection bodies of star bodies}. 

Note that the right-hand
side of (\ref{intbodyofstar}) can be written using (\ref{volume=spherradon}):
$$\|\xi\|_K^{-1}= \frac{1}{n-1} \int_{S^{n-1}\cap \xi^\bot} \|\theta\|_L^{-n+1} d\theta=
\frac{1}{n-1} R(\|\cdot\|_L^{-n+1})(\xi),$$
where $R$ is the spherical Radon transform. It means that a star body $K$ is 
the intersection body of a star body if and only if the function $\|\cdot\|_K^{-1}$
is the spherical Radon transform of a continuous positive function on $S^{n-1}.$
This allows us to introduce a more general class of bodies. We say that a star body
$K$ in $\R^n$ is an {\it intersection body} 
if there exists a finite Borel measure \index{intersection body}
$\mu$ on the sphere $S^{n-1}$ so that $\|\cdot\|_K^{-1}= R\mu$ as functionals on 
$C(S^{n-1}),$ i.e. for every continuous function $f$ on $S^{n-1},$
\begin{equation} \label{defintbody}
\int_{S^{n-1}} \|x\|_K^{-1} f(x)\ dx = \int_{S^{n-1}} Rf(x)\ d\mu(x).
\end{equation}

Intersection bodies played the crucial role in the solution of the Busemann-Petty problem
due to the following connection found by Lutwak \cite{L}: if $K$ in an origin-symmetric 
intersection body in $\R^n$ and $L$ is any origin-symmetric star body in $\R^n,$
then the inequalities $S_K(\xi)\le S_L(\xi)$ for all $\xi\in S^{n-1}$ imply that
$\vol_n(K)\le \vol_n(L),$ i.e. the answer to the Busemann-Petty problem in this
situation is affirmative. For more information about intersection bodies, see \cite[Ch.4]{K4},
\cite{KY}, \cite[Ch.8]{G3} and references there.

In this section we prove the stability of Lutwak's connection. First, we need some simple 
facts about the $\Gamma$-function.

\begin{lemma} \label{gammafunction} For any $n\in \N,$ the following inequalities hold:
$$ 1\le \frac{\left(\Gamma(\frac n2+1)\right)^{\frac{n-1}n}}{\Gamma(\frac {n+1}2 )} \le \sqrt{e},$$

$$ \frac{\Gamma(\frac{n-1}2)}{\left(\Gamma(\frac n2)\right)^{\frac{n-1}n}} \le \frac{n^{\frac{n-1}n}2^{\frac 1n}}{n-1}$$
and 
$$\sqrt{\frac n2} \le \frac{\Gamma(\frac n2+1)}{\Gamma(\frac {n+1}2)}\le \sqrt{\frac{n+1}2}.$$
\end{lemma}

\pf The the first inequality see for example \cite[Lemma 2.1]{KL}. The second inequality is a simple modification of 
the lower estimate in the first, using the property $\Gamma(x+1)=x\Gamma(x)$ of the $\Gamma$-function.

The third inequality follows from log-convexity of the $\Gamma$-function (see \cite[p.30]{K4}):
$$\Gamma^2\left(\frac n2+1\right) \le \Gamma\left(\frac {n+3}2\right)\Gamma\left(\frac {n+1}2 \right)
= \left(\frac {n+1}2 \right)\Gamma^2\left(\frac {n+1}2 \right),$$
and 
$$\Gamma^2\left(\frac {n+1}2\right) \le \Gamma\left(\frac {n}2+1\right)\Gamma\left(\frac {n}2 \right)
= \frac {2}n \Gamma^2\left(\frac {n}2+1 \right).\qed$$
\bigbreak

\begin{theorem} \label{main-int} Suppose that $\e>0$,  $K$ and $L$ are origin-symmetric
star bodies in $\R^n,$ and $K$ is an intersection body.  If for every $\xi\in S^{n-1}$
\begin{eqnarray}\label{sect1}
S_K(\xi)\le S_L(\xi) +\e,
\end{eqnarray}
then
$$\vol_n(K)^{\frac{n-1}n}  \le \vol_n(L)^{\frac{n-1}n} + \e.$$
\end{theorem}

\pf   By (\ref{volume=spherradon}), the condition (\ref{sect1}) can be written as
\begin{equation} \label{radon-buspetty}
R(\|\cdot\|_K^{-n+1})(\xi) \le R(\|\cdot\|_L^{-n+1})(\xi) + (n-1)\e,
\ \forall \xi\in S^{n-1}.
\end{equation}
Since $K$ is an intersection body, there exists a finite Borel measure 
$\mu$ on $S^{n-1}$ such that $\|\cdot\|_K^{-1}= R\mu$ as functionals on $C(S^{n-1}).$
Together with (\ref{polar-volume}), (\ref{radon-buspetty}) and the definition of $R\mu$, the latter implies that
$$n \vol_n(K) = \int_{S^{n-1}} \|x\|_K^{-1}\|x\|_K^{-n+1}\ dx $$
$$= \int_{S^{n-1}} \|x\|_K^{-n+1}\ d(R\mu)(x)
= \int_{S^{n-1}} R\left(\|\cdot\|_K^{-n+1}\right)(\xi)\ d\mu(\xi)$$
$$\le \int_{S^{n-1}} R\left(\|\cdot\|_L^{-n+1}\right)(\xi)\ d\mu(\xi) + (n-1)\e \int_{S^{n-1}}  d\mu(\xi)$$
\begin{equation} \label{eq11}
= \int_{S^{n-1}} \|x\|_K^{-1}\|x\|_L^{-n+1}\ dx + (n-1)\e \int_{S^{n-1}} d\mu(x).
\end{equation}
We estimate the first term in (\ref{eq11}) using H\"older's inequality:
$$ \int_{S^{n-1}} \|x\|_K^{-1}\|x\|_L^{-n+1}\ dx \le  \left(\int_{S^{n-1}} \|x\|_K^{-n}\ dx\right)^{\frac1n}
\left( \int_{S^{n-1}} \|x\|_L^{-n}\ dx\right)^{\frac{n-1}n}$$
\begin{equation}\label{eq12}
 = n \vol_n(K)^{\frac1n} \vol_n(L)^{\frac{n-1}n}.
 \end{equation}
 We now estimate the second term in (\ref{eq11}) adding the Radon transform of  the unit constant
 function under the integral ($R1(x)=\left|S^{n-2}\right|$ for every $x\in S^{n-1}$), 
 using again the fact that $\|\cdot\|_K^{-1}=R\mu$ and then applying H\"older's  inequality:
\begin{equation} \label{eq55}
(n-1)\e \int_{S^{n-1}} d\mu(x) = \frac{(n-1)\e}{\left|S^{n-2}\right|} \int_{S^{n-1}} R1(x)\ d\mu(x)
\end{equation}
$$=\frac{(n-1)\e}{\left| S^{n-2} \right| } \int_{S^{n-1}} \|x\|_K^{-1}\ dx $$
\begin{equation}\label{eq13}
\le  \frac{(n-1)\e}{\left|S^{n-2}\right|} \left|S^{n-1}\right|^{\frac{n-1}n} \left(\int_{S^{n-1}} \|x\|_K^{-n}\ dx\right)^{\frac1n},
\end{equation}
where 
$$\left|S^{n-2}\right|= \frac{2\pi^{\frac{n-1}2}}{\Gamma(\frac{n-1}2)} \qquad {\rm and}\qquad 
\left|S^{n-1}\right|= \frac{2\pi^{\frac{n}2}}{\Gamma(\frac{n}2)} $$ 
are the surface areas of the unit spheres in $\R^{n-1}$ and $\R^n,$ correspondingly.

We get that the quantity in (\ref{eq13}) is equal to 
$$ \frac{(n-1) \Gamma(\frac {n-1}2)}{2^{\frac 1n} \left( \Gamma(\frac n2)\right)^{\frac {n-1}n}}\e
\left(n\vol_n(K)\right)^{\frac1n} \le n\e \left(\vol_n(K)\right)^{\frac1n}$$
by the second inequality of Lemma \ref{gammafunction}.

Combining the latter inequality with (\ref{eq11}) and  (\ref{eq12}),
$$n \vol_n(K) \le n \vol_n(K)^{\frac1n} \vol_n(L)^{\frac{n-1}n} + n\e \left(\vol_n(K)\right)^{\frac1n}.\qed$$
\bigbreak
It is known that for $2\le n \le 4$ every origin symmetric convex body in $\R^n$ is an intersection body
(see \cite{G2}, \cite{Z3}, \cite{GKS} or \cite[p. 73]{K4}).  This means that the result of Theorem \ref{main-int} holds in these
dimensions for arbitrary origin-symmetric convex bodies $K,L.$ Moreover, interchanging $K,L$ in 
Theorem \ref{main-int}, we prove

\begin{co}  \label{n4} If $2\le n \le 4,$ then for any origin-symmetric convex bodies $K,L$ in $\R^n,$ 
$$\left|\vol_n(K)^{\frac{n-1}n} - \vol_n(L)^{\frac{n-1}n}\right| \le \|S_K-S_L\|_{C(S^{n-1})}.$$
\end{co}
\bigbreak

We now prove the linear separation property of Lutwak's connection.
Denote by  $$r(K) = \frac{\min_{\xi\in S^{n-1}} \rho_K(\xi)}{\vol_n(K)^{1/n}}$$
the normalized inradius of $K.$
\begin{theorem} \label{main-int1} Let $K$ and $L$ be origin-symmetric
star bodies in $\R^n$ and $\e>0.$ Assume that $K$ is an intersection body.  
If for every $\xi\in S^{n-1}$
\begin{eqnarray}\label{sect2}
S_K(\xi)\le S_L(\xi) - \e,
\end{eqnarray}
then
$$\vol_n(K)^{\frac{n-1}n}  \le \vol_n(L)^{\frac{n-1}n} -  \sqrt{\frac{2\pi}{n+1}}\ r(K) \e.$$
\end{theorem}
\pf The proof goes along the same lines as that of Theorem \ref{main-int},
with the difference that now we need a lower estimate
in place of the upper estimate (\ref{eq13}). Similarly to (\ref{eq11}) and (\ref{eq12}), we get
\begin{equation} \label{same}
n \vol_n(K) \le n \vol_n(K)^{\frac1n} \vol_n(L)^{\frac{n-1}n} - (n-1)\e \int_{S^{n-1}} d\mu(x).
\end{equation}
Similarly to (\ref{eq55}),
$$(n-1)\e \int_{S^{n-1}} d\mu(x) = \frac{(n-1)\e}{\left| S^{n-2} \right| } \int_{S^{n-1}} \|x\|_K^{-1}\ dx,$$
and, since $\|x\|_K^{-1}=\rho_K(x)$ for $x\in S^{n-1},$ using the definition of $r(K)$ we estimate the latter by
$$\ge \frac{\e (n-1)r(K) \vol_n(K)^{\frac1n} \left|S^{n-1}\right|}{\left|S^{n-2}\right|}$$
$$ = \e r(K) n \vol_n(K)^{\frac1n} \frac{(n-1)\pi^{\frac n2}\Gamma(\frac{n-1}2)}{n\pi^{\frac{n-1}2} \Gamma(\frac n2)}$$
(we multipiled and divided by $n$ and now use $\Gamma(x+1)=x\Gamma(x)$ and the third inequality 
of Lemma \ref{gammafunction})
$$= \e r(K) n \vol_n(K)^{\frac1n}\sqrt{\pi} \frac{\Gamma(\frac{n+1}2)}{\Gamma(\frac n2+1)} \ge
\e r(K) n \vol_n(K)^{\frac1n} \sqrt{\frac{2\pi}{n+1}}.$$
Combining this with (\ref{same}), we get
$$n \vol_n(K) \le n \vol_n(K)^{\frac1n} \vol_n(L)^{\frac{n-1}n} - n\e r(K) \vol_n(K)^{\frac1n} \sqrt{\frac{2\pi}{n+1}}.\qed$$
\bigbreak

We now pass to stability in the comparison result from [KYY].  The goal here is to establish
stability of volume comparison in dimensions higher than 4 without the assumption 
that $K$ is an intersection body.  We use the techniques of the Fourier approach
to sections of convex bodies that has recently been developed; see \cite{K4} and \cite{KY}. 

The Fourier transform of a
distribution $f$ is defined by $\langle\hat{f}, \phi\rangle= \langle f, \hat{\phi} \rangle$ for
every test function $\phi$ from the Schwartz space $ \mathcal{S}$ of rapidly decreasing infinitely
differentiable functions on $\R^n$. For any even distribution $f$, we have $(\hat{f})^\wedge
= (2\pi)^n f$.

If $K$ is a star body  and $0<p<n,$
then $\|\cdot\|_K^{-p}$  is a locally integrable function on $\R^n$ and represents a distribution. 
Suppose that $K$ is infinitely smooth, i.e. $\|\cdot\|_K\in C^\infty(S^{n-1})$ is an infinitely differentiable 
function on the sphere. Then by \cite[Lemma 3.16]{K4}, the Fourier transform of $\|\cdot\|_K^{-p}$  
is an extension of some function $g\in C^\infty(S^{n-1})$ to a homogeneous function of degree
$-n+p$ on $\R^n.$ When we write $\left(\|\cdot\|_K^{-p}\right)^\wedge(\xi),$ we mean $g(\xi),\ \xi \in S^{n-1}.$
If $K,L$ are infinitely smooth star bodies, the following spherical version of Parseval's
formula was proved in \cite{K5} (see \cite[Lemma 3.22]{K4}):  for any $p\in (-n,0)$
\begin{equation}\label{parseval}
\int_{S^{n-1}} \left(\|\cdot\|_K^{-p}\right)^\wedge(\xi) \left(\|\cdot\|_L^{-n+p}\right)^\wedge(\xi) =
(2\pi)^n \int_{S^{n-1}} \|x\|_K^{-p} \|x\|_L^{-n+p}\ dx.
\end{equation}

A distribution is called {\it positive definite} if its Fourier transform is a positive distribution in
the sense that $\langle \hat{f},\phi \rangle \ge 0$ for every non-negative test function $\phi.$
The following was proved in \cite{KYY}: 

\begin{lemma}\label{Lem:pos-def} {\bf (\cite[Lemma 2.3]{KYY})}
Let $K$ be an origin-symmetric convex body in $\R^n$.
Assume that $\alpha\in [n-4,n-1)$, then $\|x\|_K^{-1}\cdot |x|^{-\alpha}_2$ is a positive definite
distribution on $\R^n$.
\end{lemma}
If $K$ is infinitely smooth, by  Lemma \ref{Lem:pos-def} and \cite[Lemma 3.16]{K4}, the Fourier transform
$(|x|^{-\alpha}_2\|x\|_K^{-1})^\wedge$ is an extension of a non-negative infinitely differentiable 
function on $S^{n-1}$ to the whole $\R^n.$

\begin{theorem} \label{main} Let $\e>0,\ \alpha \in [n-4,n-1)$,  and let $K$ and $L$ be origin-symmetric
infinitely smooth convex bodies in $\R^n$, $n\ge 4$, so that for every $\xi\in S^{n-1}$
\begin{eqnarray}\label{eqn:condition11}
(-\Delta)^{\alpha/2} S_K(\xi)\le (-\Delta)^{\alpha/2} S_L(\xi) +\e.
\end{eqnarray}
Then
$$\vol_n(K)^{\frac{n-1}n}  \le \vol_n(L)^{\frac{n-1}n} + c \e,$$
where 
$$c=c(\alpha,n) = \frac{\sqrt{\pi}(n-1)\Gamma(\frac{n-\alpha-1}2)}{2^{\alpha+\frac1n} n^{\frac{n-1}n}
\Gamma(\frac{\alpha+1}2) \left(\Gamma(\frac n2)\right)^{\frac{n-1}n}}.$$
\end{theorem}

\pf  It was proved in \cite{K1} that
\begin{eqnarray}\label{eqn:defS}
S_K(\xi)=\frac{1}{\pi(n-1)}(\|x\|_K^{-n+1})^\wedge(\xi),\qquad \forall \xi\in S^{n-1}.
\end{eqnarray}
Extending $S_K(\xi)$ to $\R^n$ as a homogeneous function of degree $-1$ and using the definition
of fractional powers of the Laplacian
we get
$$
(-\Delta)^{\alpha/{2}}S_K(\theta)= \frac{1}{\pi(n-1)}(|x|^\alpha_2\|x\|_K^{-n+1})^\wedge(\theta),
$$
therefore
$$(2\pi)^n n \vol_n(K) = (2\pi)^n \int_{S^{n-1}}\|x\|_K^{-n+1}\|x\|_K^{-1}dx$$
$$=(2\pi)^n \int_{S^{n-1}}(|x|^{-\alpha}_2\|x\|_K^{-1})(|x|^{\alpha}_2\|x\|_K^{-n+1})dx$$
$$=\int_{S^{n-1}}(|x|^{-\alpha}_2\|x\|_K^{-1})^\wedge(\theta)(|x|^{\alpha}_2\|x\|_K^{-n+1})^\wedge(\theta)d\theta$$
$$={\pi(n-1)}\int_{S^{n-1}}(|x|^{-\alpha}_2\|x\|_K^{-1})^\wedge(\theta)(-\Delta)^{{\alpha}/{2}}S_K(\theta)d\theta.$$
Here we used Parseval's formula on the sphere (\ref{parseval}). 
By Lemma \ref{Lem:pos-def},
$(|x|^{-\alpha}_2\|x\|_K^{-1})^\wedge$ is a non-negative function on $S^{n-1}$, 
and we can use (\ref{eqn:condition11}) to estimate the latter quantity:
$$\le {\pi(n-1)}\int_{S^{n-1}}(|x|^{-\alpha}_2\|x\|_K^{-1})^\wedge(\theta)(-\Delta)^{{\alpha}/{2}}S_L(\theta)d\theta$$
\begin{equation} \label{eq21}
+ {\pi(n-1)}\e \int_{S^{n-1}}(|x|^{-\alpha}_2\|x\|_K^{-1})^\wedge(\theta)d\theta.
\end{equation}
Repeating the above calculation in the opposite order, we get that the first summand in (\ref{eq21}) is equal to
\begin{equation}\label{eq22}
(2\pi)^n \int_{S^{n-1}}\|x\|_L^{-n+1}\|x\|_K^{-1}dx \le (2\pi)^n n\vol_n(K)^{\frac1n} \vol_n(L)^{\frac{n-1}n}
\end{equation}
by H\"older's inequality. 

To estimate the second summand in (\ref{eq21}), we use the formula for the Fourier transform 
(in the sense of distributions; see \cite[p.194]{GS})
$$\left(|x|_2^{-n+\alpha+1}\right)^\wedge(\theta) = \frac{2^{\alpha+1}\pi^{\frac n2}\Gamma(\frac{\alpha+1}2)}
{\Gamma(\frac{n-\alpha-1}2)} |\theta|_2^{-\alpha-1}.$$
Again using Parseval's formula and then H\"older's inequality,
$$\int_{S^{n-1}}(|x|^{-\alpha}_2\|x\|_K^{-1})^\wedge(\theta)d\theta$$
$$ = 
\frac{\Gamma(\frac{n-\alpha-1}2)}{2^{\alpha+1}\pi^{\frac n2}\Gamma(\frac{\alpha+1}2)}
\int_{S^{n-1}}(|x|^{-\alpha}_2\|x\|_K^{-1})^\wedge(\theta)\left(|x|_2^{-n+\alpha+1}\right)^\wedge(\theta)d\theta$$
$$=  \frac{(2\pi)^n \Gamma(\frac{n-\alpha-1}2)}{2^{\alpha+1}\pi^{\frac n2}\Gamma(\frac{\alpha+1}2)}
\int_{S^{n-1}} \|x\|_K^{-1}\ dx$$
$$ \le \frac{(2\pi)^n \Gamma(\frac{n-\alpha-1}2)\left|S^{n-1}\right|^{\frac{n-1}n}}{2^{\alpha+1}\pi^{\frac n2}\Gamma(\frac{\alpha+1}2)}
\left(\int_{S^{n-1}} \|x\|_K^{-n}\ dx\right)^{\frac 1n}$$
$$= \frac{(2\pi)^n \Gamma(\frac{n-\alpha-1}2)\left|S^{n-1}\right|^{\frac{n-1}n}}{2^{\alpha+1}\pi^{\frac n2}\Gamma(\frac{\alpha+1}2)}\left(n\vol_n(K)\right)^{\frac 1n}$$
Combining this with (\ref{eq21}) and (\ref{eq22}), we get
$$(2\pi)^n n \vol_n(K) \le (2\pi)^n n\vol_n(K)^{\frac1n} \vol_n(L)^{\frac{n-1}n} $$
$$+ \frac{(2\pi)^n\e \pi(n-1)n^{\frac 1n}  \Gamma(\frac{n-\alpha-1}2)
\left|S^{n-1}\right|^{\frac{n-1}n}}{2^{\alpha+1}\pi^{\frac n2}\Gamma(\frac{\alpha+1}2)}\left(\vol_n(K)\right)^{\frac 1n},$$
which implies the result.
\endpf
For $\alpha<n-4$ the statement of Theorem \ref{main} is no longer true, simply because
the comparison result itself does not hold, as shown in \cite{KYY}.
\bigbreak
The corresponding separation result looks as follows:
\begin{theorem} \label{main2} Let $\e>0,\ \alpha \in [n-4,n-1)$,  $K$ and $L$ be origin-symmetric
infinitely smooth convex bodies in $\R^n$, $n\ge 4$, so that for every $\xi\in S^{n-1}$
\begin{eqnarray}\label{eqn:condition}
(-\Delta)^{\alpha/2} S_K(\xi)\le (-\Delta)^{\alpha/2} S_L(\xi) -\e.
\end{eqnarray}
Then
$$\vol_n(K)^{\frac{n-1}n}  \le \vol_n(L)^{\frac{n-1}n} - c \e,$$
where 
$$c= r(K)\frac{\pi(n-1)\Gamma(\frac{n-\alpha-1}2)}{n2^\alpha \Gamma(\frac{\alpha+1}2)\Gamma(\frac n2)}.$$
\end{theorem}

\pf Following the proof of Theorem \ref{main}, we get
$$(2\pi)^n n \vol_n(K) \le (2\pi)^n n\vol_n(K)^{\frac1n} \vol_n(L)^{\frac{n-1}n}$$
$$ - {\pi(n-1)}\e \int_{S^{n-1}}(|x|^{-\alpha}_2\|x\|_K^{-1})^\wedge(\theta)d\theta.$$
The difference with the proof of Theorem \ref{main} is that now we have to estimate
$$\int_{S^{n-1}}(|x|^{-\alpha}_2\|x\|_K^{-1})^\wedge(\theta)d\theta$$
from below. In the same way as in Theorem \ref{main} we write this integral as
$$  \frac{(2\pi)^n \Gamma(\frac{n-\alpha-1}2)}{2^{\alpha+1}\pi^{\frac n2}\Gamma(\frac{\alpha+1}2)}
\int_{S^{n-1}} \|x\|_K^{-1}\ dx.$$
The latter integral is greater or equal to $r(K) \left(\vol_n(K)\right)^{\frac 1n}\left|S^{n-1}\right|.$
The result follows.
\endpf

\section{Stability and separation for projections}

We need several more definitions from convex geometry. We refer the reader
to \cite{S2} for details.

The {\it support function} of a convex body $K$ in $\R^n$ is defined by
$$h_K(x) = \max_{\{\xi\in \R^n:\|\xi\|_K=1\}} (x,\xi),\quad x\in \R^n.$$ 
If $K$ is origin-symmetric, then $h_K$ is a norm on $\R^n.$

The {\it surface area measure} $S(K, \cdot)$ of a convex body $K$ in 
$\R^n$ is defined as follows: for every Borel set $E \subset S^{n-1},$ 
$S(K,E)$ is equal to Lebesgue measure of the part of the boundary of $K$
where normal vectors belong to $E.$ 
We usually consider bodies with absolutely continuous surface area measures.
A convex body $K$ is said to have the {\it curvature function} 
$$ f_K: S^{n-1} \to \R,$$
if its surface area measure $S(K, \cdot)$ is absolutely 
continuous with respect to Lebesgue measure $\sigma_{n-1}$ on 
$S^{n-1}$, and
$$
\frac{d S(K, \cdot)}{d \sigma_{n-1}}=f_K \in L_1(S^{n-1}),
$$
so $f_K$ is the density of $S(K,\cdot).$

By the approximation argument of \cite[Th. 3.3.1]{S2},
we may assume in the formulation of Shephard's problem that the bodies 
$K$ and $L$ are such that  their support functions $h_K,\ h_L$ are 
infinitely smooth functions on $\R^n\setminus \{0\}$.
Using \cite[Lemma 3.16]{K4}
we get in this case that
the Fourier transforms $\widehat{h_K},\ \widehat{h_L}$ are the
extensions of infinitely differentiable functions on the sphere
to homogeneous distributions on $\R^n$ of degree $-n-1.$
Moreover, by a similar approximation argument (see also \cite[Section 5]{GZ}),
we may assume that  our bodies have absolutely continuous surface area 
measures. Therefore, in the rest of this section, $K$ and $L$ are 
convex symmetric bodies with infinitely smooth support functions and absolutely 
continuous surface area measures.

The following version of Parseval's formula was proved in \cite{KRZ} (see also \cite[Lemma 8.8]{K4}):
\begin{equation} \label{pars-proj}
\int_{S^{n-1}} \widehat{h_K} (\xi) \widehat{f_L}(\xi)\ d\xi =
(2\pi)^n \int_{S^{n-1}} h_K(x) f_L(x)\ dx.
\end{equation}

The volume of a body can be expressed in terms of its support function and 
curvature function:
\begin{equation}\label{vol-proj}
\vol_n(K) = \frac 1n \int_{S^{n-1}}h_K(x) f_K(x)\ dx.
\end{equation}

If $K$ and $L$ are two convex bodies in $\R^n$ the {\it mixed volume} $V_1(K,L)$
is equal to  
$$V_1(K,L)= \frac{1}{n} \lim_{\e\to +0}
\frac{\vol_n(K+\epsilon L)-\vol_n(K)}{\e}.$$
We use the following
first Minkowski inequality (see \cite[p.23]{K4}):  
for any convex bodies $K,L$ in $\R^n,$ 
\begin{equation} \label{firstmink}
V_1(K,L) \ge \vol_n(K)^{(n-1)/n} \vol_n(L)^{1/n}.
\end{equation}
The mixed volume can also be expressed in terms of the support and
curvature functions:

\begin{equation}\label{mixvol-proj}
V_1(K,L) = \frac 1n \int_{S^{n-1}}h_L(x) f_K(x)\ dx.
\end{equation}

Let $K$ be an origin-symmetric convex body in $\R^n.$ The {\it
projection body} $\Pi K$ of $K$ is defined as an origin-symmetric convex 
body in $\R^n$ whose support function in every direction is equal to
the volume of the hyperplane projection of $K$ to this direction: 
for every $\theta\in S^{n-1},$ 
\begin{equation} \label{def:proj}
h_{\Pi K}(\theta) = \vol_{n-1}(K\vert\theta^{\perp}).
\end{equation}
If $L$ is the projection body of some convex body, we simply say 
that $L$ is a projection body. 

Both Petty \cite{Pe} and Schneider \cite{S1} in their solutions of the Shephard problem (see the introduction)
used the connection with projection bodies:  if the body $L$ (with greater projections)
is a projection body then the answer to the question of the Shephard problem is affirmative  
for any body $K.$ We now prove the stability of this connection.

Define the normalized circumradius of $L$ by
$$R(L) = \frac{\max_{\xi\in S^{n-1}} \rho_L(\xi)}{\vol_n(L)^{\frac 1n}}.$$

\begin{theorem} \label{main-proj} Suppose that $\e>0$,  $K$ and $L$ are origin-symmetric
convex bodies in $\R^n,$ and $L$ is a projection body.  If for every $\xi\in S^{n-1}$
\begin{eqnarray}\label{proj1}
P_K(\xi)\le P_L(\xi) +\e,
\end{eqnarray}
then
$$\vol_n(K)^{\frac{n-1}n}  \le \vol_n(L)^{\frac{n-1}n} + \sqrt{\frac{2\pi}n}\ R(L) \e.$$
\end{theorem}

\pf  It was proved in \cite{KRZ} that
\begin{equation} \label{f-proj}
P_K(\xi) = -\frac 1{\pi} \widehat{f_K}(\theta),\qquad \forall \xi\in S^{n-1},
\end{equation}
where $f_K$ is extended from the sphere to a homogeneous function of degree 
$-n-1$ on the whole $\R^n,$ and the Fourier transform $\widehat{f_K}$ is the 
extension of a continuous function $P_K$ on the sphere to a homogeneous of degree 1
function on $\R^n.$

Therefore, the condition (\ref{proj1}) can be written as
\begin{equation} \label{fourier-proj}
-\frac 1{\pi} \widehat{f_K}(\xi) \le  -\frac 1{\pi} \widehat{f_L}(\xi) + \e, \qquad \forall \xi\in S^{n-1}.
\end{equation}

It was also proved in \cite{KRZ} that an infinitely smooth origin-symmetric convex body 
$L$ in $\R^n$ is a projection body if and only if 
$\widehat{h_L} \le 0$ on the sphere $S^{n-1}.$ Therefore, integrating (\ref{fourier-proj})
with respect to a negative density,
$$\int_{S^{n-1}} \widehat{h_L}(\xi) \widehat{f_L}(\xi)\ d\xi \ge \int_{S^{n-1}} \widehat{h_L}(\xi) \widehat{f_K}(\xi)\ d\xi 
+ \pi\e \int_{S^{n-1}} \widehat{h_L}(\xi)\ d\xi.$$
Using this, (\ref{vol-proj}) and (\ref{pars-proj}), we get
$$ (2\pi)^n n \vol_n(L) = (2\pi)^n \int_{S^{n-1}} h_L(x) f_L(x)\  dx =
\int_{S^{n-1}} \widehat{h_L}(\xi) \widehat{f_L}(\xi)\ d\xi$$
$$\ge \int_{S^{n-1}} \widehat{h_L}(\xi) \widehat{f_K}(\xi)\ d\xi  + \pi\e\int_{S^{n-1}}\widehat{h_L}(\xi)\ d\xi$$
\begin{equation} \label{eq31}
=(2\pi)^n \int_{S^{n-1}} h_L(x) f_K(x)\  dx + \pi\e\int_{S^{n-1}}\widehat{h_L}(\xi)\ d\xi.
\end{equation}
We estimate the first summand from below using the first Minkowski inequality:
\begin{equation} \label{eq32}
(2\pi)^n \int_{S^{n-1}} h_L(x) f_K(x)\  dx \ge (2\pi)^n n \left(\vol_n(L)\right)^{\frac 1n} \left(\vol_n(K)\right)^{\frac {n-1}n}.
\end{equation}
To estimate the second summand in (\ref{eq31}), note that, by (\ref{f-proj}),
the Fourier transform of the curvature function of the Euclidean ball
$$\widehat{f_2}(\xi) = -\pi \vol_{n-1}(B_2^{n-1}) =
- \frac{\pi^{\frac{n+1}2}}{\Gamma(\frac{n+1}2)},\qquad \forall \xi\in S^{n-1},$$
where $B_2^{n-1}$ is the unit Euclidean ball in $\R^{n-1}.$
Therefore, 
$$\pi \e \int_{S^{n-1}}\widehat{h_L}(\xi)\ d\xi = - \e \frac{\Gamma(\frac{n+1}2)}{\pi^{\frac{n-1}2}}
\int_{S^{n-1}}\widehat{h_L}(\xi)\widehat{f_2}(\xi)\ d\xi$$
$$= - (2\pi)^n \e \frac{\Gamma(\frac{n+1}2)}{\pi^{\frac{n-1}2}}
\int_{S^{n-1}} h_L(x) f_2(x)\ dx$$
$$ = - (2\pi)^n \e \frac{\Gamma(\frac{n+1}2)}{\pi^{\frac{n-1}2}}
\int_{S^{n-1}} h_L(x) \ dx $$
$$\ge - (2\pi)^n \e \frac{\Gamma(\frac{n+1}2)}{\pi^{\frac{n-1}2}} R(L) \left(\vol_n(L)\right)^{\frac 1n}\left|S^{n-1}\right|$$
$$= - (2\pi)^n n \e \frac{\sqrt{\pi} R(L) \left(\vol_n(L)\right)^{\frac 1n}\Gamma(\frac{n+1}2)}{\Gamma(\frac n2 +1)},$$
where we again used Parseval's formula, the fact that $f_2=1,$ and a simple estimate 
$h_L(x)\le R(L)\left(\vol_n(L)\right)^{\frac 1n}.$

Combining this with (\ref{eq31}) and (\ref{eq32}), and using the third inequality of Lemma \ref{gammafunction},
we get
$$(2\pi)^n n \vol_n(L) \ge (2\pi)^n n \left(\vol_n(L)\right)^{\frac 1n} \left(\vol_n(K)\right)^{\frac {n-1}n} $$
$$- (2\pi)^n n \sqrt{\frac{2\pi}n} R(L) \left(\vol_n(L)\right)^{\frac 1n} \e,$$
which finishes the proof.
\endpf
We now prove the corresponding separation result.

\begin{theorem} \label{main-proj1} Suppose that $\e>0$,  $K$ and $L$ are origin-symmetric
convex bodies in $\R^n,$ and $L$ is a projection body.  If for every $\xi\in S^{n-1}$
\begin{eqnarray}\label{proj2}
P_K(\xi)\le P_L(\xi) - \e,
\end{eqnarray}
then
$$\vol_n(K)^{\frac{n-1}n}  \le \vol_n(L)^{\frac{n-1}n} - \frac{ \e}{\sqrt{e}}.$$
\end{theorem}

\pf Similarly to the proof of Theorem \ref{main-proj}, we get (\ref{eq31}), but with negative sign
in front of $\e:$
\begin{equation} \label{eq33}
(2\pi)^n n \vol_n(L) \ge (2\pi)^n \int_{S^{n-1}} h_L(x) f_K(x)\  dx - \pi\e\int_{S^{n-1}}\widehat{h_L}(\xi)\ d\xi.
\end{equation}
The difference with the proof of Theorem \ref{main-proj} is that now we need an upper estimate for
$$\pi\e\int_{S^{n-1}}\widehat{h_L}(\xi)\ d\xi = - (2\pi)^n \e \frac{\Gamma(\frac{n+1}2)}{\pi^{\frac{n-1}2}}
\int_{S^{n-1}} h_L(x) f_2(x)\ dx.$$
Using the first Minkowski inequality (\ref{firstmink}), the latter is 
$$\le  - (2\pi)^n n \e \frac{\Gamma(\frac{n+1}2)}{\pi^{\frac{n-1}2}} \left(\vol_n(L)\right)^{\frac 1n}
 \left(\vol_n(B_2^n)\right)^{\frac {n-1}n}$$
 $$= - (2\pi)^n n \e \frac{\Gamma(\frac{n+1}2)}{\left(\Gamma(\frac n2 +1)\right)^{\frac {n-1}n}} \left(\vol_n(L)\right)^{\frac 1n}
 \le - \frac{(2\pi)^n n \e}{\sqrt{e}} \left(\vol_n(L)\right)^{\frac 1n}$$
 by the first inequality of Lemma \ref{gammafunction}. In conjunction with (\ref{eq33}) and (\ref{firstmink}), (\ref{mixvol-proj}), 
 this implies the result.
 \endpf
 
 Finally, we formulate the stability version of the result from \cite{Y} mentioned in the introduction, which treats projections 
 in arbitrary dimension without the additional assumption that $L$ is a projection body. The proof follows the lines of
 the proof of Theorem \ref{main-proj} with changes corresponding to those in the proof of Theorem \ref{main}; we leave 
 this proof to the willing reader, as  well as the separation result in this case. Let us just mention that 
 one has to use the fact that for every $\alpha\in [n,n+1)$ the distribution
 $|x|_2^{-\alpha}h_L(x)$ is positive definite, which is explained in \cite{Y}. 
 
 \begin{theorem} \label{main-1} Let $\e>0,\ \alpha \in [n,n+1)$,  $K$ and $L$ be origin-symmetric
infinitely smooth convex bodies in $\R^n$, $n\ge 3$, so that for every $\xi\in S^{n-1}$
$$(-\Delta)^{\alpha/2} P_L(\xi)\le (-\Delta)^{\alpha/2} P_K(\xi) +\e.$$
Then
$$\vol_n(K)^{\frac{n-1}n}  \le \vol_n(L)^{\frac{n-1}n} + c \e,$$
where 
$$c=\frac{\Gamma(\frac{n-\alpha+1}2) \left|S^{n-1}\right| R(L)}{2^{\alpha+1}\pi^{\frac n2}
\Gamma(\frac{\alpha+1}2)n}.$$
\end{theorem}
Note that this is no longer true if $\alpha<n,$ because the underlying comparison result fails, 
as shown in \cite{Y}.
\bigbreak
{\bf Acknowledgement.} The author wishes to thank
the US National Science Foundation for support through 
grants DMS-0652571 and DMS-1001234.

\end{document}